\documentclass[12pt,a4paper]{article}

\usepackage{amsthm}
\usepackage{thmtools}

\declaretheorem[style=plain,qed=$\square$,numberwithin=section]{theorem}

\declaretheorem[style=plain,qed=$\square$,sibling=theorem]{corollary}
\declaretheorem[style=plain,qed=$\square$,sibling=theorem]{proposition}
\declaretheorem[style=plain,qed=$\square$,sibling=theorem]{claim}

\declaretheorem[style=definition,qed=$\circ$,sibling=theorem]{assumption}
\declaretheorem[style=definition,qed=$\circ$,sibling=theorem]{definition}
\declaretheorem[style=definition,qed=$\circ$,sibling=theorem]{remark}

\declaretheoremstyle[
  spaceabove=6pt, spacebelow=6pt,
  headfont=\normalfont\slshape,
  notefont=\normalfont\slshape, notebraces={(}{)},
  bodyfont=\normalfont,
  postheadspace=1em,
  qed=$\blacksquare$
]{proofstyle}

\usepackage{graphicx}
\usepackage{amssymb}
\usepackage{epstopdf}

\usepackage{xcolor}
\usepackage{graphicx}

\usepackage{amsfonts,amssymb,amsmath}
\usepackage{mathtools}
\usepackage{mathrsfs}
\usepackage{dsfont}
\usepackage{centernot}
\DeclareMathAlphabet{\mathpzc}{OT1}{pzc}{m}{it}

\usepackage[english]{babel}
\usepackage[utf8]{inputenc}
\usepackage{hyperref}

\hypersetup{
colorlinks=true, linktocpage=true, pdfstartpage=3, pdfstartview=FitV,
breaklinks=true, pdfpagemode=UseNone, pageanchor=true, pdfpagemode=UseOutlines,
plainpages=false, bookmarksnumbered, bookmarksopen=true, bookmarksopenlevel=1,
hypertexnames=true, pdfhighlight=/O, urlcolor=purple, linkcolor=black, citecolor=green!50!black,
pdftitle={},
pdfauthor={D.D. Baptista de Souza, H. Stein Shiromoto},
pdfsubject={},
pdfkeywords={},
pdfcreator={pdfLaTeX},
pdfproducer={TeXShop}
}

\usepackage{longtable}

\usepackage{todonotes}

\usepackage{enumerate}
\usepackage{paralist}

\usepackage{setspace} 
\usepackage{multirow}

\usepackage{caption}
\input{commands}

\newcommand{\bbN}{\mathbb{N}}

\newcommand{\R}{\mathbb{R}}
\newcommand{\Rs}{\mathbb{R}_{>0}}
\newcommand{\Ras}{\mathbb{R}_{\geq0}}

\newcommand{\bbZ}{\mathbb{Z}}

\newcommand{\bB}{\mathbf{B}}

\newcommand{\bD}{\mathbf{D}}

\newcommand{\bM}{\mathbf{M}}
\newcommand{\bN}{\mathbf{N}}

\newcommand{\bR}{\mathbf{R}}
\newcommand{\bS}{\mathbf{S}}

\newcommand{\bY}{\mathbf{Y}}

\newcommand{\bZ}{\mathbf{Z}}

\newcommand{\cC}{\mathcal{C}}

\newcommand{\cK}{\mathcal{K}}
\newcommand{\cL}{\mathcal{L}}

\newcommand{\cP}{\mathcal{P}}

\newcommand{\dz}{\dot{z}}

\usepackage{authblk}

\begin{document}

\title{A Region-Dependent Gain Condition for Asymptotic Stability}
\author[1]{Humberto Stein Shiromoto\thanks{The work of the first and the third authors are partly supported by HYCON2 Network of Excellence Highly-Complex and Networked Control Systems, grant agreement 257462. E-mail: \texttt{humberto.shiromoto@ieee.org}.}}
\author[2,3]{Vincent Andrieu}
\author[4]{Christophe Prieur}
\affil[1]{\small INCAS3, Dr. Nassaulaan 9, 9401 HJ Assen, The Netherlands}
\affil[2]{\small Fachbereich C - Mathematik und Naturwissenschaften, Bergische Universit\"at Wuppertal,
Gaußstraße 20, 42097 Wuppertal, Germany}
\affil[3]{\small LAGEP, Universit\'e de Lyon, Rue Victor Grignard, CPE, B\^at. G, 69622 Villeurbanne, France.}
\affil[4]{\small GIPSA-lab, Grenoble Campus, 11 rue des Math\'ematiques, BP 46,
38402 Saint Martin d'H\`eres Cedex, France}

\date{}

\maketitle
\thispagestyle{empty}
\pagestyle{empty}


\begin{abstract}
	A sufficient condition for the stability of a system resulting from the interconnection of dynamical systems is given by the small gain theorem. Roughly speaking, to apply this theorem, it is required that the gains composition is continuous, increasing and upper bounded by the identity function. In this work, an alternative sufficient condition is presented for the case in which this criterion fails due to either lack of continuity or the bound of the composed gain is larger than the identity function. More precisely, the local (resp. non-local) asymptotic stability of the origin (resp.  global attractivity  of a compact set) is ensured by a  region-dependent  small gain condition. Under an additional condition  that implies convergence of solutions for almost all initial conditions  in a suitable domain, the  almost  global asymptotic stability of the origin is ensured. Two examples illustrate and motivate this approach.
\end{abstract}


\section{Introduction}

	The use of nonlinear input-output gains for stability analysis was introduced in \cite{Zames1966P1} by considering a system as an input-output operator. The condition that ensures stability, called Small Gain Theorem, of interconnected systems is based on the contraction principle.

The work \cite{Sontag1989} introduces a new concept of gain relating the input to system states. This notion of stability links Zames' and Lyapunov's approaches \cite{Sontag:2001}. Characterizations in terms of dissipation and Lyapunov functions are given in \cite{SontagWang1995}.

In \cite{Jiangetal1994}, the contraction principle is used in the input-to-state stability notion to obtain an equivalent Small Gain Theorem. A formulation of this criterion in terms of Lyapunov functions may be found in \cite{Jiangetal:1996}.

Besides stability analysis, the Small Gain Theorem may also be used for the design of dynamic feedback laws satisfying robustness constraints. The interested reader is invited to see \cite{FreemanKokotovic2008,Sastry:1999} and references therein. Other versions of the Small Gain theorem do exist in the literature, see \cite{Angeli:2007,AstolfiPraly:2012,Ito2006,ItoJiang:2009} for not necessarily ISS systems.

	In order to apply the Small Gain Theorem, it is required that the composition of the nonlinear gains is smaller than the argument for all of its positive values. Such a condition, called Small Gain Condition, restricts the application of the Small Gain Theorem to a composition of well chosen gains.

In this work, an alternative criterion for the stabilization of interconnected systems is provided when a single Small Gain Condition does not hold globally. It consists in showing that if the two conditions hold: 1) a local (resp. non-local) Small Gain Condition holds in a local (resp. non-local) region of the state space, and the intersection of the local and non-local is empty, and 2) outside the union of these regions, the set of initial conditions from which the associated trajectories do not converge to the local region has measure zero, then the resulting interconnected system is almost asymptotically stable (this notion is precisely defined below). In this paper, a sufficient condition guaranteeing this property to hold is presented.
Moreover, for planar systems, an extension of the Bendixson's criterion to regions which are not simply connected is given. This allows to obtain global asymptotic stability of the origin.

This approach may be seen as a unification of two small gain conditions that hold in different regions: a local and a non-local. The use of a unifying approach for local and non-local properties is well known in the literature see \cite{AndrieuPrieur2010} in the context of control Lyapunov functions, see \cite{Chaillet:2012} when uniting iISS and ISS properties.

This paper is organized as follows. In Section \ref{sec:introduction:Preliminaires}, the system under consideration and the problem statement are presented. Section \ref{sec:asymmetric:assumptions} states the assumptions to solve the problem under consideration and the main results. Section \ref{sec:asymmetric:illustration} presents examples that illustrate the assumptions and main results. In Section \ref{sec:asymmetric:proof} the proofs of the main results are presented. Section \ref{sec:conclusion} collects some concluding remarks.

{\small\textbf{Notation.} Let $k\in\bbZ_{>0}$. Let $\bS$ be a subset of $\R^k$ containing the origin, the notation $\bS_{\neq0}$ stands for $\bS\setminus\{0\}$. The \emph{closure} of $\bS$ is denoted by  $\cl\{\bS\}$. Let $x\in\R^k$, the notation $|x|$ stands for \emph{Euclidean norm} of $x$.  An \emph{open} (resp. \emph{closed}) \emph{ball centered at} $x\in\R^k$ \emph{with radius} $r>0$ is denoted by $\bB_{<r}(x)$ (resp. $\bB_{\leq r}(x)$). A continuous function $f:\R^k\to\R$ is \emph{positive definite} if, for every $x\in\R^k\setminus\{0\}$, $f(x)>0$ and $f(0)=0$.  It is \emph{proper} if $|f(x)|\to\infty$, as $|x|\to\infty$. By $\lloc^\infty(\R,\R^k)$ the class of functions $\eta:\R\to\R^k$ that are \emph{locally essentially bounded}. By $\cC^s$ it is denoted the class of \emph{$s$-times continuously differentiable functions}, by $\cP$ it is denoted the class of \emph{positive definite functions}, by $\cK$ it is denoted the class of \emph{continuous, positive definite  and strictly increasing functions} $\gamma:\Ras\to\Ras$; it is denoted by $\cK_\infty$ if, in addition, they are \emph{unbounded}. Let $c\in\Rs$, the notation $\Omega_{\diamond c}(f)$ stands for the \emph{subset of} $\R^k$ defined by $\{x\in\R^k:f(x)\diamond c\}$, where $\diamond$ is a \emph{comparison operator} (i.e., $=$, $<$, $\geq$ etc). The \emph{support} of the function $f$ is the set $\supp:=\{x\in\R^k:f(x)\neq0\}$. By $\lloc^\infty(\Ras,\R^k)$ it is denoted the class of functions $g:\Ras\to\R^k$ that are locally essentially bounded. Let $x,\bar{x}\in\Ras$, the notation $x\nearrow\bar{x}$ (resp. $x\searrow\bar{x}$) stands for $x\to\bar{x}$ with $x<\bar{x}$ (resp. $x>\bar{x}$).}

\section{Background and problem statement}\label{sec:introduction:Preliminaires}

Consider the system
\begin{equation}\label{eq:x subsystem}
			\dx(t)=f(x(t),u(t)),
\end{equation}
where, for every $t\in\Ras$, $x(t)\in\R^n$, and $u\in\lloc^\infty(\Ras,\R^m)$, for some positive integers $n$ and $m$. Also, $f\in\cC^1(\R^{n+m},\R^n)$.  A solution of \eqref{eq:x subsystem} with initial condition $x$, and input $u$ at time $t$ is denoted by $X(t,x,u)$. From now on, arguments $t$ will be omitted, and assume that the origin is \emph{input-to-stable stable} (ISS for short) \emph{for \eqref{eq:x subsystem}}. For further details on this concept, the interested reader is invited to consult \cite{Sontag:2001} or \cite{SontagWang1996}.

A locally Lipschitz function $V:\R^n\to\Ras$ for which there exist $\underline{\alpha}_x$, $\overline{\alpha}_x\in\cK_\infty$ such that, for every $x\in\R^n$, $\underline{\alpha}_x(|x|)\leq V(x)\leq \overline{\alpha}_x(|x|)$ is called \emph{storage function}.

Inspired by \cite{Dashkovskiy:2010a,Liberzon:2013}, the following notion of derivative will be used.

\begin{definition}
	Consider the function $\xi:[a,b)\to\R$,  the limit at $t\in[a,b)$ 
	\begin{equation*}
		D^+\xi(t)=\limsup_{\tau\searrow0}\tfrac{\xi(t+\tau)-\xi(t)}{\tau}
	\end{equation*}
	(if it exists) is called \emph{Dini derivative}. Let $k_1$ and $k_2$ be positive integers, $(y_1,y_2)\in\R^{k_1}\times\R^{k_2}$, functions $\varphi:\R^{k_1+k_2}\to\R$, $h_1:\R^{k_1}\to\R^{k_1}$ and $h_2:\R^{k_2}\to\R^{k_2}$. The limit
	\begin{equation*}
		D^+_{h_1,h_2}\varphi(y_1,y_2)=\limsup_{\tau\searrow0}\tfrac{\varphi(y_1 + \tau h_1(y_1),y_2 + \tau h_2(y_2))-\varphi(y_1,y_2)}{\tau}
	\end{equation*}
	(if it exists) is called \emph{Dini derivative of $\varphi$ in the $h_1$ and $h_2$-directions at $(y_1,y_2)$}.\footnote{When the Dini derivative is taken in only one direction, the subscript denotes only such a direction.}
\end{definition}

 If, for a given storage function $V$, there exist  a proper function $\lambda_x\in(\cC^0\cap\cP)(\R^n,\Ras)$, and $\alpha_x\in\cK_{\infty}$  called \emph{ISS-Lyapunov gain}  such that, for every $(x,u)\in\R^n\times\R^m$,
\begin{equation}\label{eq:dini Lyap inequality:V}
	|x|\geq\alpha_x(|u|)\Rightarrow D_f^+V(x,u)\leq-\lambda_x(x),
\end{equation}
then $V$ is called \emph{ISS-Lyapunov function} for \eqref{eq:x subsystem}. As in \cite{Dashkovskiy:2010a}, the proof that the existence of an ISS-Lyapunov implies that \eqref{eq:x subsystem} is ISS goes along the lines presented in \cite{SontagWang1995}. 

Consider the system\footnote{A solution of \eqref{eq:z subsystem} with initial condition $z$, and input $v$ at time $t$ is denoted by $Z(t,z,v)$.}
\begin{equation}\label{eq:z subsystem}
			\dot{z}=g(v,z),
\end{equation}
where $v\in\lloc^\infty(\Ras,\R^n)$, $z\in\R^m$, and $g\in\cC^1(\R^{n+m},\R^m)$. From now on, assume that $W:\R^{n+m}\to\Ras$ is an ISS-Lyapunov function for \eqref{eq:z subsystem} with $\lambda_z\in(\cC^0\cap\cP)(\R^m,\Ras)$, and $\alpha_z\in\cK_{\infty}$ satisfying, for every $(v,z)\in\R^{n+m}$, 
\begin{equation}\label{eq:dini Lyap inequality:W}
	W(z)\geq\alpha_z(|v|)\Rightarrow D_g^+W(v,z)\leq-\lambda_z(z).
\end{equation}\black

{\bfseries System under consideration.} Interconnecting systems \eqref{eq:x subsystem} and \eqref{eq:z subsystem} yields the system
\begin{equation}\label{eq:general system}
	\left\{\begin{array}{rcl}
			\dx&=&f(x,z),\\
			\dz&=&g(x,z).
	\end{array}\right.
\end{equation}
 Using the vectorial notation $y=(x,z)$, system \eqref{eq:general system} is denoted by $\dot{y}=h(y)$. A solution initiated from $y$ in $\R^{n+m}$ and evaluated at time $t$ is denoted $Y(t,y)$. The two ISS-Lyapunov inequalities \eqref{eq:dini Lyap inequality:V} and \eqref{eq:dini Lyap inequality:W} can be rephrased as follows. For every couple $(x,z)\in\R^{n+m}$,
\begin{equation}\label{eq:ISS implications}
	\begin{array}{rclcrcl}
		V(x)&\geq&\gamma(W(z))&\Rightarrow&D_f^+V(x,z)&\leq&-\lambda_x(x),\\
		W(z)&\geq&\delta(V(x))&\Rightarrow&D_g^+W(x,z)&\leq&-\lambda_z(z)
	\end{array}
\end{equation}
with suitable functions $\gamma,\delta\in\cK_{\infty}$.

A sufficient condition that ensures the stability of \eqref{eq:general system} is given by the small gain theorem \cite{Jiangetal:1996}. Roughly speaking if,
\begin{equation}\label{eq:introduction:general sgc}
	\forall s\in\Rs,\quad\gamma\circ\delta(s)<s,
\end{equation}
then the origin is globally asymptotically stable for \eqref{eq:general system}.

{\bfseries Problem statement.} At this point, it is possible to explain the problem under consideration. ISS systems for which \eqref{eq:introduction:general sgc} does not hold in a bounded set of $\Ras$ are considered. This paper shows that by merging small gain arguments in different regions of the state space and employing some tools from measure theory, a sufficient condition ensuring almost global asymptotic stability of the origin is possible to be given. For planar interconnected systems, by using an extension of Bendixon's criterion, global asymptotic stability of the origin may be established. 

\section{Assumptions and main results}\label{sec:asymmetric:assumptions}

\begin{assumption}\label{hyp:Lyapunov inequalities} There exist constant values $0\leq\underline{M}<\overline{M}\leq\infty$ and $0\leq\underline{N}<\overline{N}\leq\infty$, and class $\mathcal{K}_{\infty}$ functions $\gamma$ and $\delta$ such that, for every $(x,z)\in\bS\subset\R^n\times\R^m$, the implications
\begin{eqnarray}
			V(x)\geq \gamma(W(z))&\Rightarrow& D_f^+V(x,z)\leq-\lambda_x(x)\label{eq:Lyap ineq:x}\\
			W(z)\geq \delta(V(x))&\Rightarrow& D_g^+W(x,z)\leq-\lambda_z(z)\label{eq:Lyap ineq:z}
\end{eqnarray}
hold, where
\begin{equation}\label{eq:S}
	\begin{array}{rcl}
\bS&:=&\{(x,z)\in\R^n\times\R^m:\underline{M}\leq V(x)\leq \overline{M},\\
&&W(z)\leq\overline{N}\}
\cup\{(x,z)\in\R^n\times\R^m:\\
&& V(x)\leq\overline{M},
\underline{N}\leq W(z)\leq \overline{N}\},
	\end{array}
\end{equation}
\end{assumption}

In other words, Assumption \ref{hyp:Lyapunov inequalities} states that the set $\Omega_{\leq\underline{M}}(V)\times\Omega_{\leq\underline{N}}(W)$ is \emph{locally} ISS for the $x$ and $z$-subsystems of \eqref{eq:general system}. To see more details on locally ISS systems, the interested reader may consult \cite{Dashkovskiy:2010}.

\begin{assumption}\label{hyp:sgc}\hfill
\begin{equation}\label{eq:sgc:interval}
	\begin{array}{rcl}
		\text{if}\ \overline{M}<\infty,&s\in[\underline{M},\overline{M}]\setminus\{0\},&\gamma\circ\delta(s)<s,\\
	\text{if}\ \overline{M}=\infty,&s\in[\underline{M},\overline{M})\setminus\{0\},&\gamma\circ\delta(s)<s.
	\end{array}
\end{equation}
\end{assumption}

Assumption \ref{hyp:sgc} states that the small gain condition holds in the interval corresponding  to the value of $V$, when $x$ is restricted to $\bS$.

\begin{proposition}\label{prop:attractivity in intervals}
	 Under Assumptions \ref{hyp:Lyapunov inequalities} and \ref{hyp:sgc}, if
\begin{equation}\label{eq:max < min}
		\widetilde{M}:=\max\{\gamma^{-1}(\underline{M}),\underline{N}\}<\min\{\delta(\overline{M}),\overline{N}\}=:\widehat{M},
\end{equation}
	 then there exists a proper function $U\in\cP(\R^{n+m},\Ras)$ that is locally Lipschitz on $\R^{n+m}\setminus\{0\}$ and such that,
\begin{equation*}
	\forall y\in\Omega_{\leq\widehat{M}}(U)\setminus\Omega_{\leq\widetilde{M}}(U),\quad\displaystyle\limsup\limits_{t\to\infty}   U(Y(t,y))\leq\widetilde{M}.
\end{equation*}
Moreover, if $\gamma,\delta\in(\cC^1\cap\cK_\infty)$, then a suitable $U$ can be defined, for every $(x,z)\in\R^n\times\R^m$, by
\begin{equation}\label{eq:proposition:function U}
	U(x,z)=\max\left\{\tfrac{\delta(V(x))+\gamma^{-1}(V(x))}{2},W(z)\right\}.
\end{equation}
\end{proposition}

Condition \eqref{eq:max < min} implies that $\Omega_{\leq\widetilde{M}}(U)\subsetneq\Omega_{\leq\widehat{M}}(U)$. Proposition \ref{prop:attractivity in intervals} states that solutions of \eqref{eq:general system} starting in $\Omega_{\leq\widehat{M}}(U)$ will converge to the set $\Omega_{\leq\widetilde{M}}(U)$. The proof of Proposition \ref{prop:attractivity in intervals} is provided in Section \ref{sec:proof of prop:attractivity in intervals}.

\begin{corollary}{[Local stabilization]}\label{cor:local stability}  
	Consider Assumptions \ref{hyp:Lyapunov inequalities} and \ref{hyp:sgc} with the constant values $\underline{M}=\underline{N}=0$, $M_\ell:=\overline{M}<\infty$ or $N_\ell:=\overline{N}<\infty$. The set $\Omega_{\leq\widehat{M}_\ell}(U_\ell)$ is included in the basin of attraction of the origin of \eqref{eq:general system}, where $U_\ell$ and $\widehat{M}_\ell$ are given by Proposition \ref{prop:attractivity in intervals}.
\end{corollary}

 In other words, Corollary \ref{cor:local stability} states that the set $\Omega_{\leq\widehat{M}_\ell}(U_\ell)$ is an estimation of the set of initial conditions from which issuing solution of \eqref{eq:general system} remain close and converge to the origin.

 Before stating the second corollary, some concepts regarding the asymptotic behaviour of solutions are recalled.  A set $\bM\subset\R^{n+m}$ is said to be \emph{positively  invariant} with respect to \eqref{eq:general system} if, for every $t\in\Ras$, $y\in\bM\Rightarrow Y(t,y)\in\bM$ (cf. \cite[p. 127]{Khalil:2001}).  A compact positively invariant set $\bM\subset\R^{n+m}$ is said to be \emph{globally attractive} if, for all $y\in\R^{n+m}$, $\lim_{t\to\infty}|Y(t,y)|_\bM=0$.

\begin{corollary}{[Global attractivity]}\label{cor:global attractivity} 
	Consider Assumptions \ref{hyp:Lyapunov inequalities} and \ref{hyp:sgc} with the constant values $M_g:=\underline{M}>0$ or $N_g:=\underline{N}>0$, and $\overline{M}=\overline{N}=\infty$. The set $\Omega_{\leq\widetilde{M}_g}(U_g)$ is globally attractive for \eqref{eq:general system}, where $U_g$ and $\widetilde{M}_g$ are given by Proposition \ref{prop:attractivity in intervals}.
\end{corollary}

 In other words, Corollary \ref{cor:global attractivity} states that the set $\Omega_{\leq\widetilde{M}_g}(U_g)$ is an estimation of the global attractor of \eqref{eq:general system}.

 The proofs of Corollaries \ref{cor:local stability} and \ref{cor:global attractivity} are not provided and follow from Proposition \ref{prop:attractivity in intervals}. The interested reader may also consult \cite{Dashkovskiy:2010a,Dashkovskiy:2010}.

\begin{figure}
	\centering
	\scalebox{0.75}{\input{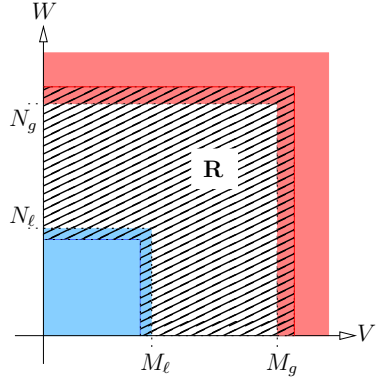}}
	\caption{Illustration of sets $\Omega_{\leq M_\ell}(V)\times\Omega_{\leq N_\ell}(W)$ (blue region), $\Omega_{=\widetilde{M}_\ell}(U_\ell)$ (dark blue line), $\Omega_{\geq\underline{M}}(V)\times\Omega_{\geq\underline{N}}(W)$ (pink region), $\Omega_{=\widetilde{M}_g}(U_g)$ (red line), and $\bR=\cl\{\Omega_{\leq\widehat{M}_g}(U_g)\setminus\Omega_{\leq\widetilde{M}_\ell}(U_\ell)\}$ (pattern filled).}
	\label{fig:R}
\end{figure}

 Under the assumptions of Corollaries \ref{cor:local stability} and \ref{cor:global attractivity}, if the estimation of the global attractor $\Omega_{\leq\widetilde{M}_g}(U_g)$ is contained in the estimation of the basin of attraction $\Omega_{\leq\widehat{M}_\ell}(U_\ell)$, then global asymptotic stability of the origin for \eqref{eq:general system} follows trivially. However, when this inclusion does not hold, the set $\bR=\Omega_{\leq\widetilde{M}_g}(U_g)\setminus\Omega_{\leq\widehat{M}_\ell}(U_\ell)$ is not empty, and solutions of \eqref{eq:general system} may converge to positively invariant sets contained in $\bR$ instead (cf. Birkhoff's Theorem \cite{Isidori:1995}). Figure \ref{fig:R} illustrates the region $\bR$ obtained in this situation.

The next result provides sufficient conditions ensuring that, for almost every initial condition, issuing solutions remain close and converge to the origin. For the case in which \eqref{eq:general system} is planar, global asymptotic stability of the origin is established.

Before stating the main results,  the concept of stability introduced in \cite{Angeli:2004} is presented.  The origin is called \emph{almost globally asymptotically stable} for \eqref{eq:general system} if it is locally stable in the Lyapunov sense and attractive for almost every initial condition. More precisely, there exists $\aleph\subset\R^{n+m}$, with $\mu(\aleph)=0$ such that, for every $y\in\R^{n+m}\setminus\aleph$, $\lim_{t\to\infty}|Y(t,y)|=0$, where $\mu$ is the Lebesgue measure.

\begin{theorem}\label{thm:almost global stability}
	 Under Assumptions \ref{hyp:Lyapunov inequalities} and \ref{hyp:sgc}, if the constant values of Corollaries \ref{cor:local stability} and \ref{cor:global attractivity} are such that $M_\ell<M_g$ or $N_\ell<N_g$, there exists a function $\rho\in\cC^1(\R^{n+m}\setminus\{0\},\Ras)$ with $\supp(\rho)\supseteq\bR$, where $\bR=\cl\{\Omega_{\leq \widetilde{M}_g}(U_g)\setminus\Omega_{\leq\widehat{M}_\ell}(U_\ell)\}$, and if for every $y\in\bR$, $\div (h\rho)(y)>0$, then the origin is almost globally asymptotically stable for \eqref{eq:general system}.
\end{theorem}

 In other words, Theorem 1 states that with an extra assumption on the vector field of system (4), solutions converge to the origin for almost every initial condition and the origin is locally asymptotically stable.  The proof of Theorem \ref{thm:almost global stability} is provided in Section \ref{sec:proof of thm:almost global stability}.

\begin{theorem}\label{thm:global stability}
	 Let $n=m=1$. Under Assumptions \ref{hyp:Lyapunov inequalities} and \ref{hyp:sgc}, if the constant values of Corollaries \ref{cor:local stability} and \ref{cor:global attractivity} are such that $M_\ell<M_g$ or $N_\ell<N_g$, and for every $y\in\bR=\cl\{\Omega_{\leq \widetilde{M}_g}(U_g)\setminus\Omega_{\leq\widehat{M}_\ell}(U_\ell)\}$, $\diverg h(y)\neq0$ and $h(y)\neq0$, then the origin is globally asymptotically stable for \eqref{eq:general system}.
\end{theorem}

 Theorem \ref{thm:global stability} states that, when \eqref{eq:general system} is planar and under mild conditions on the vector field, the origin is globally asymptotically stable for \eqref{eq:general system}. In other words, no $\omega$-limit sets exist in $\bR$.  The proof of Theorem \ref{thm:global stability} is provided in Section \ref{sec:proof of thm:global stability}.

\section{Illustration}\label{sec:asymmetric:illustration}

 Here, the results given in the previous section are illustrated in two examples. The first concerns the vectorial case, while the second concerns the planar case.

\subsection{A class of systems satisfying an asymptotic small-gain condition}

Recall system \eqref{eq:general system}, and assume that there exist locally Lipschitz and proper functions $V\in \cP(\R^n,\Ras)$ and $W\in\cP(\R^m,\Ras)$ satisfying\footnote{From Remark 2.4 of \cite{SontagWang1995} Eq. \eqref{eq:system of Lyap ineqs} is equivalent to \eqref{eq:ISS implications}.}, for every $(x,z)\in\R^n\times\R^m$,
\begin{equation}\label{eq:system of Lyap ineqs}
	\left\{\begin{array}{rcccl}
				D^+_fV(x,z)&\leq&-V(x)&+&\gamma(W(z)),\\
				D^+_gW(x,z)&\leq&-W(z)&+&\delta(V(x)),
			\end{array}\right.
\end{equation}
where $\gamma,\delta\in\cK_\infty$ are such that\footnote{Note that, in contrast to \cite{AstolfiPraly:2012,Ito2006,ItoJiang:2009}, no information is given about the behaviour of the function $s\mapsto\gamma\circ\delta(s)$, in the interval $(0,\infty)$.}
\begin{equation}\label{eq:asymptotic small-gain}
	\lim_{s\to0}\frac{\gamma\circ\delta(s)}{s}<1\quad\text{and}\quad\lim_{s\to\infty}\frac{\gamma\circ\delta(s)}{s}<1.
\end{equation}
In other words, the composition $\gamma\circ\delta$ satisfies the small-gain condition in the bi-limit: $0$ and $\infty$. Note that to apply \cite{AstolfiPraly:2012} it would be necessary to impose that, for every $b_1,b_2\in\Ras$, and for some $\varepsilon>0$, $\gamma(b_1)\delta(b_2)\leq (1-\varepsilon)b_1b_2$.

From \eqref{eq:asymptotic small-gain}, there exists a positive constant $M_\ell$ (resp. $M_g$) that is sufficiently small (resp. large) and such that, for every $s\in(0,M_\ell]$ (resp. $s\in[M_g,\infty)$), $\gamma\circ\delta(s)<s$. Together with \eqref{eq:system of Lyap ineqs} and since $W$ is continuous and proper, Assumptions \ref{hyp:Lyapunov inequalities} and \ref{hyp:sgc} hold locally on the compact set $\Omega_{\leq M_\ell}(V)\times\Omega_{\leq N_\ell}(W)$ (resp. non-locally on the set $\Omega_{\geq M_g}(V)\times\Omega_{\geq N_g}(W)$). 

Since condition Eq. \eqref{eq:max < min} is satisfied, as formulated in Corollary \ref{cor:local stability} (resp. \ref{cor:global attractivity}), from this result the set $\Omega_{\leq\widehat{M}_\ell}(U_\ell)\subset\R^{n+m}$ (resp. $\Omega_{\leq\widetilde{M}_g}(U_g)\subset\R^{n+m}$) is an estimation of the basin of attraction of the origin (resp. global attractor) of \eqref{eq:example:planar:general system}. Also, $\Omega_{\leq\widehat{M}_\ell}(U_\ell)\subset\Omega_{\leq\widetilde{M}_g}(U_g)$.

Let\footnote{Note that, in contrast to \cite{Angeli:2007}, here it is not assumed that $\Omega_{\geq\widetilde{M}_g}(U_g)\cup\Omega_{\leq\widehat{M}_\ell}(U_\ell)=\R^{n+m}$.}
\begin{equation*}
	\left\{\begin{array}{rclcl}
				r_V(x,z)&:=&-V(x)&+&\gamma(W(z)),\\
				r_W(x,z)&:=&-W(z)&+&\delta(V(x)).
			\end{array}\right.
\end{equation*}
For any $(x,z)$ belonging to the sets $\Omega_{\leq\widehat{M}_\ell}(U_\ell)$ and $\Omega_{\geq\widetilde{M}_g}(U_g)\subset\R^{n+m}$ either \\
\begin{inparaenum}[1.]
\item $r_V(x,z)\geq0$ and $r_W(x,z)\leq0$ or;\\
\item $r_V(x,z)\leq0$ and $r_W(x,z)\geq0$ or;\\
\item $r_V(x,z)\leq0$ and $r_W(x,z)\leq0$.\\
\end{inparaenum}
While in the compact set $\bR=\Omega_{\leq\widetilde{M}_g}(U_g)\setminus\Omega_{\leq\widehat{M}_\ell}(U_\ell)$, it may happen that $r_V(x,z)\geq0$ and $r_V(x,z)\geq0$. In this case, the result given in \cite{Angeli:2007} can not be applied, because it requires that the union of the regions described by items 1-3 forms a cover of $\R^{n+m}$.

Note that in contrast to \cite{Ito2006}, the existence of a Lyapunov function candidate for the system \eqref{eq:general system} whose derivative is definite negative on $\R^{n+m}$ is not requested.

Let, for every $(x,y)\in\R^{n+m}$, $\rho(x,y)=(V(x)+W(z))^{-1}$, and assume that for every $(x,z)\in\bR$,
\begin{equation}\label{eq:example:div condition}
(\div h(x,z))(V(x)+W(z))\geq\gamma(W(z))+\delta(V(x)).
\end{equation}
In such a compact set, note that
\begin{equation*}
	\begin{array}{rcl}
		\div(h\rho)(x,z)&=&\rho(x,z)\div h(x,z)+\grad\rho(x,z)\cdot h(x,z)\\
		&=&\dfrac{\div h(x,z)}{V(x)+W(z)}-\dfrac{D^+_fV(x,z)+D^+_gW(x,z)}{(V(x)+W(z))^2}\\
		&\geq&\dfrac{\div h(x,z)+1}{V(x)+W(z)}-\dfrac{\gamma(W(z))+\delta(V(x))}{(V(x)+W(z))^2}\\
		&>&0,
	\end{array}
\end{equation*}
where the first inequality is due to \eqref{eq:system of Lyap ineqs}: for every $(x,z)\in\R^n\times\R^m$,
\begin{equation*}
	\begin{array}{rcccl}
		-D^+_fV(x,z)-D^+_gW(x,z)&\geq& V(x)&-&\gamma(W(z))\\
		&&+W(z)&-&\delta(V(x)),
	\end{array}
\end{equation*}
and last inequality is due to \eqref{eq:example:div condition}.

From Theorem \ref{thm:almost global stability}, the origin is almost globally asymptotically stable for the system \eqref{eq:general system}, as stated below.

\begin{proposition}
	 Assume that there exist ISS-Lyapunov functions for the sub-systems of \eqref{eq:general system} with ISS-Lyapunov gains satisfying \eqref{eq:asymptotic small-gain},  assume moreover that \eqref{eq:example:div condition} holds, then the origin is almost globally asymptotically stable for system \eqref{eq:general system}.
\end{proposition}

\subsection{The planar case}

Consider the system
\begin{equation}\label{eq:example:planar:general system}
		\left\{\begin{array}{rclcccl}
			\dx&=&f(x,z)&=&-1.5x&+&2p(z),\\
			\dz&=&g(x,z)&=&-z&+&\sin(\pfrac{x^2}{10}),
		\end{array}\right.
	\end{equation}
	where, for every $s\in\R$, $p(s)\coloneqq \pfrac{s^3}{3}-3\pfrac{s^2}{2}+2s$.
	
	Let, for every $x\in\R$ (resp. $z\in\R$), $V(x)=|x|$ (resp. $W(z)=|z|$). Taking its Lie derivative in the $f$-direction yields, for every $(x,z)\in\R^2$,
	\begin{equation}\label{eq:example:planar:D_fV}
		D^+_fV(x,z)\leq -1.5V(x)+2|p(W(z))|.
	\end{equation}
	Define, for every $s\in\Ras$, $\gamma(s)=\max\{1.3|p(r)|:0\leq r\leq s\}$. From \eqref{eq:example:planar:D_fV},
	\begin{equation}\label{eq:example:planar:ISS Lyap impl:V}
		V(x)\geq\gamma(W(z))\Rightarrow D^+_fV(x,z)\leq-\lambda_x(x)
	\end{equation}
	holds with a suitable $\lambda_x\in(\cC^0\cap\cP)(\R,\Ras)$. The Lie derivative of $W$ in the $g$-direction yields, for every $(x,z)\in\R^2$,
	\begin{equation*}
		D_g^+W(x,z)\leq-W(z)+\left|\sin\left(\tfrac{V(x)^2}{10}\right)\right|.
	\end{equation*}
	which can be rephrased as follows
	\begin{equation*}\label{eq:example:planar:ISS Lyap impl:W}
		W(z)\geq\delta(V(x))\Rightarrow D_g^+W(x,z)\leq-\lambda_z(z)
	\end{equation*}
	with a suitable $\lambda_z\in(\cC^0\cap\cP)(\R,\Ras)$, where $\delta(s)=\max\{|\sin(\pfrac{r^2}{10})|:0\leq r\leq s\}$.
	
	The composition of the function $\gamma$ and $\delta$ yields
	\begin{equation*}
			\gamma\circ\delta(s)=\max \left\{ 1.3|p(r)|:0\leq r\leq \max\limits_{0\leq a\leq s}\left\{\left|\sin\left(\tfrac{a^2}{10}\right)\right|\right\}\right\}.
	\end{equation*}
	Note that there exist values $\bar{s}>0$ for which $\gamma\circ\delta(\bar{s})=1.11$. Also,
	\begin{equation}\label{eq:example:planar:asymptotic small gain condition}
		\lim\limits_{s\to0}\tfrac{\gamma\circ\delta(s)}{s}<1,\quad\text{and}\quad\lim\limits_{s\to\infty}\tfrac{\gamma\circ\delta(s)}{s}<1.
	\end{equation}
	
	From \eqref{eq:example:planar:asymptotic small gain condition}, and following the reasoning of the previous example, there exists $M_\ell>0$ small (resp. $M_g>0$ large) enough such that, for every $s\in(0,M_\ell]$ (resp. $s\in[M_g,\infty)$), $\gamma\circ\delta(s)<s$. Also, there exist\footnote{Recall that $\gamma$ and $\delta$ are continuous positive definite functions. Thus, they are strictly increasing in a neighbourhood of the origin. Note also that $\gamma$ is proper.} $\gamma_\ell,\delta_\ell\in\cK_\infty$ (resp.\footnote{Although $\gamma_g$ is of class $\cK$, the result of Proposition \ref{prop:attractivity in intervals} is still applicable. The main difference in this case would be the construction of the function $\sigma\in\cC^1\cap\cK_\infty$ satisfying \eqref{eq:proof:set glob att:delta < sigma < gamma inv}. The interested reader may consult \cite{Jiangetal:1996} to check how this is done.} $\gamma_g\in\cK$ and $\delta_g\in\cK_\infty$) satisfying, for every $s\in[0,M_\ell]$ (resp. $s\in[M_g,\infty)$), $\gamma_\ell(s)=\gamma(s)$ and $\gamma_\ell(s)=\gamma(s)$ (resp. $\gamma_g(s)=\gamma(s)$ and $\delta_g(s)=\delta(s)$). Thus, analogously to the reasoning of the previous example, Assumptions \ref{hyp:Lyapunov inequalities} and \ref{hyp:sgc} hold locally on the compact set $\Omega_{\leq M_\ell}(V)\times\Omega_{\leq N_\ell}(W)$ (resp. non-locally on the set $\Omega_{\geq M_g}(V)\times\Omega_{\geq N_g}(W)$). 
	
	Since condition Eq. \eqref{eq:max < min} is satisfied, as formulated in Corollary \ref{cor:local stability} (resp. \ref{cor:global attractivity}), from this result the set $\Omega_{\leq\widehat{M}_\ell}(U_\ell)\subset\R^{n+m}$ (resp. $\Omega_{\leq\widetilde{M}_g}(U_g)\subset\R^{n+m}$) is estimation of the basin of attraction of the origin (resp. global attractor) of \eqref{eq:example:planar:general system}. Also, $\Omega_{\leq\widehat{M}_\ell}(U_\ell)\subset\Omega_{\leq\widetilde{M}_g}(U_g)$.
	
	It now remains to check whether there exist $\omega$-limit sets in $\bR=\Omega_{\leq\widetilde{M}_g}(U_g)\setminus\Omega_{\leq\widehat{M}_\ell}(U_\ell)$. Since
	\begin{equation*}
		\begin{split}
		\tfrac{\partial f}{\partial x}(x,z)+\tfrac{\partial g}{\partial z}(x,z)\equiv-2.5\ \text{and}\\ f(x,z)=0=g(x,z)&\Leftrightarrow (x,z)=(0,0),
				\end{split}
	\end{equation*}
	from Theorem \ref{thm:global stability} the origin is globally asymptotically stable for \eqref{eq:example:planar:general system}.\black
	
\section{Proofs}\label{sec:asymmetric:proof}

\subsection{Proof of Proposition \ref{prop:attractivity in intervals}}\label{sec:proof of prop:attractivity in intervals}

\begin{proof} The proof of Proposition \ref{prop:attractivity in intervals} is based on the proof of \cite[Theorem 3.1]{Jiangetal:1996}. Here, it is divided into 3 parts. Firstly, the function $\sigma\in\cK_\infty\cap\cC^1$ is obtained. In the second part, the Dini derivative of a locally Lipschitz and proper function $U\in\cP(\R^{n+m},\Ras)$ is shown to be decreasing in the set $\bS$ defined in \eqref{eq:S}. In the third part, solutions of \eqref{eq:general system} starting in $\Omega_{\leq\widehat{M}}(U)\setminus\Omega_{\leq\widetilde{M}}(U)$ are shown to converge to $\Omega_{\leq\widetilde{M}}(U)$.

	{\scshape First Part.} Under Assumptions \ref{hyp:Lyapunov inequalities} and \ref{hyp:sgc}, the  function $\gamma$ being of class $\mathcal{K}_\infty$ satisfies, for every $s\in\Rs$, $\delta(s)<\gamma^{-1}(s)$. Together with the fact that  $\delta$ is of class $\cK_\infty$, from \cite[Lemma A.1]{Jiangetal:1996}, there exists $\sigma\in\cK_\infty\cap\cC^1$ whose derivative is strictly positive and satisfies,
\begin{equation}\label{eq:proof:set glob att:delta < sigma < gamma inv}
	\forall s\in\Rs,\quad\delta(s)<\sigma(s)<\gamma^{-1}(s).
\end{equation}
	{\scshape Second Part.} Define, for every $(x,z)\in\R^n\times\R^m$, $U(x,z)=\max\{\sigma(V(x)),W(z)\}$. Note that $U\in(\cC^0\cap\cP)(\R^{n+m},\Ras)$ is a proper function\black. Pick $({x},{z})\in\R^n\times\R^m$, one of three cases is possible: $\sigma(V({x}))<W({z})$, $W({z})<\sigma(V({x}))$ or $W({z})=\sigma(V({x}))$. The proof follows by showing that the Dini derivative of $U$ is negative definite. For each case, assume that $(x,z)\in\bS_{\neq0}:=\bS\setminus\{(0,0)\}$,  where $\bS$ is defined in \eqref{eq:S}.
	
	\emph{Case 1.} Assume that $\sigma(V(x))<W(z)$. This implies that $ U(x,z)=W(z)$ and $D^+_{f,g} U(x,z)=D_g^+W(x,z)$. From \eqref{eq:proof:set glob att:delta < sigma < gamma inv}, $\delta(V(x))<\sigma(V(x))<W(z)$. Since $(x,z)\in\bS_{\neq0}$, the inequality $D_g^+W(x,z)\leq-\lambda_z(z)$ follows from \eqref{eq:Lyap ineq:z}. Thus, $W(z)>\sigma(V(x))\Rightarrow D^+_{f,g} U(x,z)\leq -\lambda_z(z)$.
	
	\emph{Case 2.} Assume that $W(z)<\sigma(V(x))$. This implies that $ U(x,z)=\sigma(V(x))\ \text{and}\ D^+_{f,g} U(x,z)=\sigma'(V(x))D_f^+V(x,z)$. Since $(x,z)\in\bS_{\neq0}$, and from \eqref{eq:proof:set glob att:delta < sigma < gamma inv},
	\begin{equation}\label{eq:proof:set glob att:W < delta<t gamma<gamma inv}
					W(z)<\sigma(V(x))<\gamma^{-1}(V(x)).
			\end{equation}
	From \eqref{eq:Lyap ineq:x}, the inequality $D_f^+V(x,z)\leq-\lambda_x(x)$ holds.
	
	{\itshape Case 3.} Assume that $W(z)=\sigma(V(x))$. Let $U^\ast(x,z):=W(z)=\sigma(V(x))$. This implies
			\begin{equation*}
				\begin{split}
					D^{+}_{f,g} U^{\ast}(x,z)=\displaystyle\limsup_{t\searrow0}\dfrac{1}{t}[\max\{\sigma(V(X(t,x,z))),\\
					W(Z(t,z,x))\}- U^{\ast}(x,z)]\\
					=\displaystyle\limsup_{t\searrow0}\max\left\{\dfrac{\sigma(V(X(t,x,z)))-\sigma(V(x))}{t},\dfrac{W(Z(t,z,x))-W(z)}{t}\right\}\\
					=\max\{\sigma'(V(x))D_f^+V(x,z),D_g^+W(x,z)\}.
									\end{split}
			\end{equation*}
			The analysis of $D^+_{f,g} U^\ast$ is divided in two sub cases. In the first one, the function $D_g^+W$ is analyzed while in the last the function $D_f^+V$ is analyzed.
			
			\hspace{2pt}\underline{Case 3.a}. {\itshape The analysis of $D_g^+W$.} From \eqref{eq:proof:set glob att:delta < sigma < gamma inv},  and the fact that $x\neq0$ and $z\neq0$, the inequality $\delta(V(x))<\sigma(V(x))=W(z)$ holds. Moreover since $(x,z)\in\bS_{\neq0}$, the inequality $D_g^+W(x,z)\leq-\lambda_z(z)$ follows from \eqref{eq:Lyap ineq:z}.
			
			\hspace{2pt}\underline{Case 3.b}. {\itshape The analysis of $D_f^+V$.} From \eqref{eq:proof:set glob att:delta < sigma < gamma inv},  and the fact that $x\neq0$ and $z\neq0$, the inequality $W(z)=\sigma(V(x))<\gamma^{-1}(V(x))$ holds. Moreover, since $(x,z)\in\bS_{\neq0}$, the inequality $D_f^+V(x,z)\leq-\lambda_x(x)$ follows from \eqref{eq:Lyap ineq:x}.

			Summing up Case 3, $0\neq W(z)=\sigma(V(x))\Rightarrow D^{+}_{f,g} U^{\ast}(x,z)\leq-\min\{\sigma'(V(x))\lambda_x(x),\lambda_z(z)\}$.
		
		\begin{claim}\label{claim:existence of sub level set and inclusion}
			There exists $c>0$ such that $\Omega_{\leq c}(U)\subset\Omega_{\leq\overline{M}}(V)\times\Omega_{\leq\overline{N}}(W)$. Moreover, the constants $\widetilde{M}$ and $\widehat{M}$ are such that
			\begin{equation}\label{eq:inclusion}
				\begin{split}
				(\Omega_{\leq\underline{M}}(V)\times\Omega_{\leq\underline{N}}(W))\subset\Omega_{\leq\widetilde{M}}(U)\subset\Omega_{\leq\widehat{M}}(U)\\
				\subset(\Omega_{\leq\overline{M}}(V)\times\Omega_{\leq\overline{N}}(W)).
				\end{split}
			\end{equation}
		\end{claim}
		The proof of Claim \ref{claim:existence of sub level set and inclusion} is provided in Section \ref{sec:proof of claims}.
		
		From the above case study and \eqref{eq:inclusion},
	\begin{equation}\label{eq:proof:DU<-E}
		\widetilde{M}\leq U(x,z)\leq\widehat{M}\Rightarrow D^{+}_{f,g} U(x,z)\leq - E(x,z),
	\end{equation}
	where $E\in(\cC^0\cap\cP)(\R^{n+m},\R)$ is the proper function defined, for every $(x,z)\in\R^n\times\R^m$, by $E(x,z)=\min\{\sigma'(V(x))\lambda_x(x),\lambda_z(z)\}$.
	
	\textsc{Third part.} The local Lipschitz property of  $U$ on $\R^n\times\R^m\setminus\{(0,0)\}$ is due to the fact that $\sigma(V(\cdot))$ (resp. $W(\cdot)$) is locally Lipschitz on $\R^n\setminus\{0\}$ (resp. $\R^m$).
	
	 From \cite[Theorem 4.3]{Rouche:1977} and \eqref{eq:proof:DU<-E}, for all $(x,z)\in\R^n\times\R^m$, and all $t\in\Ras$, along solutions of \eqref{eq:general system}, $D^+ U(X(t,x,z),Z(t,z,x))=D^+_{f,g} U(X(t,x,z),Z(t,z,x))$.
			
			Since solutions of \eqref{eq:general system} are absolutely continuous functions and the righthand side of $E$  is a continuous and positive definite function, from \cite[Remark 4.4.b]{Rouche:1977}, for every $(x,z)$ such that $\widetilde{M}\leq U(x,z)\leq\widehat{M}$, and all $t\in\Ras$, the function
			\begin{equation}\label{eq:proof:set glob att:U<beta}
				t\mapsto  U(X(t,x,z),Z(t,z,x))
			\end{equation}
			is strictly decreasing and satisfies 
			$$U^\infty:=\lim_{t\to\infty} U(X(t,x,z),Z(t,z,x))\leq\widetilde{M}.$$
			To see this claim suppose, for purposes of contradiction, that $ U^\infty>\widetilde{M}$. From the continuity of $U$, there exists $\varepsilon>0$ such that $ U^\infty-\varepsilon>\widetilde{M}$ and $ U^\infty-\varepsilon\leq U(x,z)\leq  U^\infty+\varepsilon$. Since $U$  is proper, the constant $\xi=\min\{ E(x,z)>0:U^\infty-\varepsilon\leq U(x,z)\leq U^\infty+\varepsilon\}$ exists. Recalling the definition of $U$, there exists $T>0$ such that, for all $t\geq T$, $ U(X(t,x,z),Z(t,z,x))- U^\infty<\varepsilon$. Moreover, from the definition of the constant $\xi$,
			\begin{equation*}
				\begin{split}
					 U(X(t,x,z),Z(t,z,x))- U(X(T,x,z),Z(T,z,x))=\\
					\displaystyle\int_{T}^t D^+ U(X(s,x,z),Z(s,z,x))\,ds\leq -\xi(t-T).
				\end{split}
			\end{equation*}
			Then, 
			\begin{equation*}
				\begin{split}
					 U^\infty=\lim_{t\to\infty} U(X(t,x,z),Z(t,z,x))\\
					= U(X(T,x,z),Z(T,z,x))\\
					+\lim_{t\to\infty}\displaystyle\int_{T}^t D^+ U(X(s,x,z),Z(s,z,x))\,ds\leq-\infty
				\end{split}
			\end{equation*}
		which contradicts the fact that $U$  is positive definite. Thus, $ U^\infty\leq\widetilde{M}$. Hence, solutions of \eqref{eq:general system} starting in $\Omega_{\leq\widehat{M}}(U)\setminus\Omega_{\leq\widetilde{M}}(U)$  converge towards $\Omega_{\leq\widetilde{M}}(U)$.
		
		To see that $U$ can be given by \eqref{eq:proposition:function U}, note that $U$ relies on the computation of $\sigma$. Let, for every $s\in\Ras$, $\sigma(s)=\pfrac{(\delta(s)+\gamma^{-1}(s))}{2}$. Its derivative yields, for every $s>0$, $2\sigma'(s)=\delta'(s)+\pfrac{1}{(\gamma'\circ\gamma(s))}$ which is positive, because\footnote{Recall that $\delta,\gamma\in(\cC^1\cap\cK_\infty)$.} $\delta'(s)>0$ and $\gamma'\circ\gamma^{-1}(s)>0$. Moreover, such a function $\sigma$ satisfies \eqref{eq:proof:set glob att:delta < sigma < gamma inv}. This concludes the proof.
\end{proof}

\subsection{Proof of Theorem \ref{thm:almost global stability}}\label{sec:proof of thm:almost global stability}

 This proof is divided into 4 parts. The first one shows that solutions starting in $\Omega_{\geq\widetilde{M}_g}(U_g)$ converge to $\bR$. The second part shows that almost all solutions starting in $\bR$ converges to $\Omega_{\leq\widehat{M}_\ell}(U_\ell)$. The third part shows that solutions starting in the latter set converge to the origin. The fourth part concludes the almost global asymptotic stability of the origin.

\textsc{1st part.} From Corollary \ref{cor:global attractivity}, the set $\Omega_{\leq\widetilde{M}_g}(U_g)$ is globally attractive for \eqref{eq:general system}, where $\widetilde{M}_g=\max\{\gamma^{-1}_g(M_g),N_g\}$, $M_g$ and $N_g$ are defined in Corollary \ref{cor:global attractivity}, and $\gamma_g$ is given by Assumption \ref{hyp:Lyapunov inequalities}.

\textsc{2nd part.} From the proof of Proposition \ref{prop:attractivity in intervals}, there exist proper functions $U_g,E_g\in(\cC^0\cap\cP)(\R^{n+m},\Ras)$ (resp. $U_\ell,E_\ell\in(\cC^0\cap\cP)(\R^{n+m},\Ras)$) with $U_g$ (resp. $U_\ell$) being also locally Lipschitz and such that, for every $y\in\Omega_{\geq\widetilde{M}_g}(U_g)$, $D^+_hU_g(y)\leq-E_g(y)$ (resp. for every $y\in\Omega_{\leq\widehat{M}_\ell}(U_\ell)$, $D^+_hU_\ell(y)\leq-E_\ell(y)$). 

To see that $\Omega_{\leq\widehat{M}_\ell}(U_\ell)\subsetneq\Omega_{\leq\widetilde{M}_g}(U_g)$. From the proof of Claim \ref{claim:existence of sub level set and inclusion}, $U_\ell(x,z)\leq\widehat{M}_\ell\Rightarrow\max\{V(x),W(z)\}\leq\min\{M_\ell,N_\ell\}$. Analogously, $U_g(x,z)\geq\widetilde{M}_g$ $\Rightarrow$\linebreak  $\min\{V(x),W(z)\}\geq\max\{M_g,N_g\}$. Since $\min\{M_\ell,N_\ell\}<\max\{M_g,N_g\}$, $\Omega_{\leq\widehat{M}_\ell}(U_\ell)\subsetneq\Omega_{\leq\widetilde{M}_g}(U_g)$.

The proof proceeds by showing that, for almost every initial condition staring in $\bR=\Omega_{\leq\widetilde{M}_g}(U_g)\setminus\Omega_{\leq\widehat{M}_\ell}(U_\ell)$, issuing solutions of \eqref{eq:general system} converge to $\Omega_{\leq\widehat{M}_\ell}(U_\ell)$. To do so, the same lines as in \cite[Theorem 1]{Rantzer:2001} and \cite[Theorem 3]{Angeli:2004} are followed. However, here a less conservative condition is required, since a set that is only positively invariant, and the divergence to be positive only in a compact set are needed.
	
	Let $\bZ\subset\R^n$ a set given by\footnote{Note that $\bZ$ is the set of all initial conditions belonging to $\Omega_{\leq\widetilde{M}_g}(U_g)$ from which issuing solutions do not converge to $\Omega_{\leq\widehat{M}_\ell}(U_\ell)$.}
\begin{equation*}
	\bZ=\bigcap_{l=1}^\infty\{y\in\Omega_{\leq\widetilde{M}_g}(U_g):U_\ell(Y(t,y))>\widehat{M}_\ell,t>l\}.
\end{equation*}

 For every $t\in\R$, let $Y(t,\bZ)=\{Y(t,z):z\in\bZ, t\in \dom(z)\}$, where $\dom(z)$ is the maximum time interval where $Y(t,z)$ exists. Since $\Omega_{\leq \widetilde{M}_g}(U_g)$ is positively invariant, $\bZ$ is also positively invariant. Thus, given a fixed $\tau\in\Rs$, for all $t\geq \tau$, $Y(t,\bZ)\subset Y(\tau,\bZ)$. Hence, for all $t\in\Ras$,
\begin{equation}\label{eq:measure inequality}
	\textstyle\int_{Y(t,\bZ)}\rho(y)\,dy-\int_{\bZ}\rho(y)\,dy\leq0,
\end{equation}
where $\rho\in\cC^1(\R^{n+m}\setminus\{0\},\Ras)$ and $\supp(\rho)\supseteq\bR$.

From Liouville's Theorem (see \cite[Lemma A.1]{Rantzer:2001}), for every $t\in\Ras$,
\begin{equation*}
	\textstyle\int_0^t\int_{Y(s,\bZ)}\diverg(h\rho)(y)\,dyds=\int_{Y(t,\bZ)}\rho(y)\,dy-\int_{\bZ}\rho(y)\,dy.
\end{equation*}
Since $\bZ\subset\bR$, for every $t\in\Ras$, the inequality
\begin{equation*}
	\begin{array}{rcl}
		t\int\limits_{Y(t,\bZ)}\diverg(h\rho)(y)\,dy&\leq&\int_0^t\int_{Y(s,\bZ)}\diverg(h\rho)(y)\,dyds\\
		&\leq&\int_{Y(t,\bZ)}\rho(y)\,dy-\int_{\bZ}\rho(y)\,dy
	\end{array}
\end{equation*}
holds. From \eqref{eq:measure inequality}, for every $t\in\Ras$, $\int_{Y(t,\bZ)}\diverg(h\rho)(y)\,dy\leq0$. Together with the fact that, for every $y\in\bR$, $\diverg(h\rho)(y)>0$, it yields $\int_{Y(t,\bZ)}\diverg(h\rho)(y)\,dy=0$, for every $t\in\Ras$. Thus, for every $t\in\Ras$, $Y(t,\bZ)$ has Lebesgue measure zero. In particular, $\bZ$ has also Lebesgue measure zero. Consequently, for almost every $y\in\bR$, $\limsup_{t\to\infty}U_\ell(Y(t,y))\leq \widehat{M}_\ell$.

It remains to check if the initial conditions belonging to $\Omega_{\geq\widetilde{M}_g}(U_g)$ from which issuing solutions converge to $\bZ$ have also measure zero. Since $\bZ$ is positively invariant, for all $t_1<t_2\leq0$, $Y(t_2,\bZ)\subset Y(t_1,\bZ)$. This inclusion implies that
$
\bY:=\cup_{t\leq 0}\{Y(t,\bZ)\} = \cup_{l\in\bbZ_{<0}}\{Y(l,\bZ)\}
$.
Hence, the set $\bY$ is a countable union of images of $\bZ$ by the flow. Since $\bZ$ is measurable and, for every $t\in\dom(y)$, the map $\bZ\ni y\mapsto Y(t,y)$ is a diffeomorphism\footnote{Because \eqref{eq:general system} is of class $\cC^1$ and solutions are unique.}, $\bY$ is also measurable.

For every $t\in\dom(\bZ)$, $\int_{Y(t,\bZ)}dz \leq \int_\bZ$ $|\grad Y(t,y)|\, dy=0$,
because $\bZ$ has measure zero. This implies that, for all $t\in\dom(\bZ)$, the set $Y(t,\bZ)$ has measure zero. Since $\bY$ is the countable union of sets of measure zero, it has also measure zero.\footnote{Recall that $\bY$ is the set of initial conditions from which issuing solutions of \eqref{eq:general system} converge $\bZ$.} Hence the set of solutions starting in $\Omega_{\geq\widetilde{M}_g}(U_g)$ that converge to $\bZ$ have also measure zero.

\textsc{3rd part.} From Corollary \ref{cor:local stability}, the set $\Omega_{\leq\widehat{M}_\ell}(U_\ell)$ is contained in the basin of attraction of the origin, where $\widehat{M}_\ell=\min\{\delta_\ell(M_\ell),N_\ell\}$, $M_\ell$ and $N_\ell$ are defined in Corollary \ref{cor:local stability}, and $\gamma_\ell$ is given by Assumption \ref{hyp:Lyapunov inequalities}.

\textsc{4th part.} From the above discussion, the origin is locally stable and almost globally attractive for \eqref{eq:general system}. Thus, it is almost globally asymptotically stable for \eqref{eq:general system}. This concludes the proof.\hfill$\blacksquare$

\subsection{Proof of Theorem \ref{thm:global stability}}\label{sec:proof of thm:global stability}

Before proving Theorem \ref{thm:global stability}, some concepts regarding the asymptotic behavior of solutions of planar systems are recalled. A point $p$ is said to be a \emph{positive limit point of} $Y(\cdot,y)$ if there exists a sequence $\{t_n\}_{n\in\bbN}$, with $t_n\to\infty$ as $n\to\infty$, such that $Y(t_n,y)\to p$ as $n\to\infty$ (cf. \cite[p. 127]{Khalil:2001}). The set $\omega(y)$ of all positive limit points of $Y(\cdot,y)$ is called \emph{$\omega$-limit set of $y$} (cf. \cite[p. 517]{Isidori:1995}). For planar systems, a closed curve $C\subset\R^2$ is called \emph{closed orbit} if $C$ is not an equilibrium point and there exists a time $T<\infty$ such that, for each $y\in C$, $Y(nT,y)=y$, $\forall n\in\bbZ$ (cf. \cite[Definition 2.6]{Sastry:1999}).

\begin{proof}
	The proof of Theorem \ref{thm:global stability} follows the same line as the proof of Theorem \ref{thm:almost global stability}. The difference here consists in the second and fourth parts.  

	\textsc{1st part.} Recall that from Corollary \ref{cor:global attractivity}, the set $\Omega_{\leq\widetilde{M}_g}(U_g)$ is globally attractive for \eqref{eq:general system}, where $\widetilde{M}_g=\max\{\gamma^{-1}_g(M_g),N_g\}$, $M_g$ and $N_g$ are defined in Corollary \ref{cor:global attractivity}, and $\gamma_g$ is given by Assumption \ref{hyp:Lyapunov inequalities}.

\textsc{2nd part} (Bendixson's criterion for non simply connected regions). From the proof of Proposition \ref{prop:attractivity in intervals}, there exist proper functions $U_g,E_g\in(\cC^0\cap\cP)(\R^2,\Ras)$ (resp. $U_\ell,E_\ell\in(\cC^0\cap\cP)(\R^2,\Ras)$) with $U_g$ (resp. $U_\ell$) being also locally Lipschitz and such that, for every $y\in\Omega_{\geq\widetilde{M}_g}(U_g)$, $D^+_hU_g(y)\leq-E_g(y)$ (resp. for every $y\in\Omega_{\leq\widehat{M}_\ell}(U_\ell)$, $D^+_hU_\ell(y)\leq-E_\ell(y)$). 

 Since the set $\bR=\cl\{\Omega_{\leq \widetilde{M}_g}(U_g)\setminus\Omega_{\leq\widehat{M}_\ell}(U_\ell)\}$ is compact, and for each $y\in\bR$, $U_g(y)\neq0$, from \cite[Theorem 2.5]{Alberti:2013}\\
 \begin{inparaitem}
 	\item The set $\Omega_{=\widetilde{M}_g}(U_g)$ has finite perimeter;\\
	\item The function $U_g$ is almost every where differentiable on $\Omega_{=\widetilde{M}_g}(U_g)$;\\
	\item Let $\bN_g\subset\Omega_{=\widetilde{M}_g}(U_g)$ be set in which $U_g$ is not differentiable. There exists a Lipschitz parametrization $p_g:[a_g,b_g]\subset\R\to\Omega_{=\widetilde{M}_g}(U_g)$ that is injective and satisfies, for almost every $s\in[a_g,b_g]$, $p_g(s)\notin\bN_g$ and $p_g'(s)$ is perpendicular to $\nabla U_g(p_g(s))$.
 \end{inparaitem}
 
 Recall that by assumption, for every $y\in\bR$, $h(y)\neq0$. Together with the fact that $h\in\cC^1(\R^2)$, and almost each sublevel set of $U_g$ has finite perimeter. From the generalized divergence theorem \cite[Theorem 1.7]{Marzocchi:2005} (see also \cite{Pfeffer:2012})
	\begin{equation}\label{eq:generalized divergence theorem}
			\textstyle\iint\limits_{\Omega_{\leq \widetilde{M}_g}(U_g)}\hspace{-12.5pt} \diverg h(y)\,dy
			=\hspace{-15pt}\textstyle\oint\limits_{\Omega_{=\widetilde{M}_g}(U_g)}\hspace{-15pt} h(y)\cdot n_g(y)\,dxdz.
	\end{equation}
	Together with the above discussions and the existence of the parametrization $p_g$, for almost every $s\in[a_g,b_g]$, $h(p_g(s))\cdot n_g (p_g(s))<0$, where for almost every $s\in[a_g,b_g]$, $n_g(p_g(s))=\pfrac{\nabla U_g(p_g(s))}{|\nabla U_g(p_g(s))|}$,
	\begin{equation}\label{eq:divergence:exterior}
			\hspace{-3.5pt}\textstyle\iint\limits_{\Omega_{\leq \widetilde{M}_g}(U_g)}\hspace{-12.5pt} \diverg h(y)\,dy\hspace{-2.5pt}=\hspace{-12.5pt}\textstyle\int\limits_{[a_g,b_g]}\hspace{-6pt} h(p_g(s))\hspace{-1.5pt}\cdot\hspace{-1.5pt} n_g(p_g(s))\,ds<0.
	\end{equation}
	Analogously to the above, and by letting $p_\ell:[a_\ell,b_\ell]\to\Omega_{=\widehat{M}_\ell}(U_\ell)$ be a parametrization of $\Omega_{=\widehat{M}_\ell}(U_\ell)$ with outward unit normal $n_\ell$, based on Equation \eqref{eq:generalized divergence theorem},
	\begin{equation}\label{eq:divergence:interior}
			\textstyle\iint\limits_{\Omega_{\leq \widehat{M}_\ell}(U_\ell)}\hspace{-10pt}\diverg h(y)\,dy=\hspace{-2.5pt}\textstyle\int\limits_{[a_\ell,b_\ell]}\hspace{-7.5pt}h(p_\ell(s))\cdot n_\ell(p_\ell(s))\,ds<0.
	\end{equation}
	
	Suppose, for purposes of contradiction, that there exists a closed orbit $C\in\R^2$, parametrized by $p:[a,b]\to C$ and with outward unit normal $n$, and contained in $\bR$. From the generalized divergence theorem,
	\begin{equation}\label{eq:divergence:closed trajectory}
		\textstyle\iint\limits_{\bD_C}\hspace{-2.5pt} \diverg h(x,z)\,dxdz
		=\hspace{-5pt}
		\textstyle\int\limits_{[a,b]}\hspace{-5pt} h(p(s))\cdot n(p(s))\,ds
		=0,
	\end{equation}
	where $\bD_C$ is the simply connected region bounded by $C$.
	
	Note that,
	\begin{equation*}
		\begin{array}{rcl}
			\textstyle\iint\limits_{\Omega_{\leq\widetilde{M}_g}(U_g)\setminus\bD_C} \hspace{-22.5pt}\diverg h(y)\,dy&=&\hspace{-12.5pt}\textstyle\iint\limits_{\Omega_{\leq \widetilde{M}_g}(U_g)}\hspace{-12.5pt} \diverg h(y)\,dy
			-\textstyle\iint\limits_{\bD_C}  \diverg h(y)\,dy\\
			&=&\textstyle\iint\limits_{\Omega_{\leq \widetilde{M}_g}(U_g)}  \diverg h(y)\,dy,
		\end{array}
	\end{equation*}
	where the last equality is due to \eqref{eq:divergence:closed trajectory}. From \eqref{eq:divergence:exterior},
	\begin{equation}\label{eq:divergence:D D_C_0}
		\textstyle\iint\limits_{\Omega_{\leq\widetilde{M}_g}(U_g)\setminus\bD_C} \diverg h(y)\,dy< 0.
	\end{equation}
	On the other hand,
	\begin{equation*}
		\begin{array}{rcl}
			\textstyle\iint\limits_{\bD_C\setminus\Omega_{\leq\widehat{M}_\ell}(U_\ell)}\hspace{-20pt}\diverg h(y)\,dy&=&\textstyle\iint\limits_{\bD_C} \diverg h(y)\,dy
			-\hspace{-12.5pt}\textstyle\iint\limits_{\Omega_{\leq\widehat{M}_\ell}(U_\ell)}\hspace{-12.5pt}\diverg h(y)\,dy\\
			&=&-\textstyle\iint\limits_{\Omega_{\leq\widehat{M}_\ell}(U_\ell)} \diverg h(y)\,dy,
		\end{array}
	\end{equation*}
	where the last equality is due to \eqref{eq:divergence:closed trajectory}. From \eqref{eq:divergence:interior},
	\begin{equation}\label{eq:divergence:D_C_0 D_0}
			\textstyle\iint\limits_{\bD_C\setminus\Omega_{\leq\widehat{M}_\ell}(U_\ell)}\diverg h(x,z)\,dxdz> 0.
	\end{equation}
	
	From \eqref{eq:divergence:D D_C_0}, \eqref{eq:divergence:D_C_0 D_0} and the continuity of $\diverg h$, the function $\diverg h$ changes sign in $\bR$. Thus, there exists $\bar{y}\in\bR$ such that $\diverg h(\bar{y})=0$ which is a contradiction with the hypothesis $\div(y)\neq0$, for every $y\in\bR$. Thus, there exist no closed orbits $C$ contained in $\bR$.
	
	From the Poincar\'e-Bendixson Theorem \cite[Theorem 2.15]{Sastry:1999}, the $\omega$-limit set of a solution starting in $\bR$ is a closed orbit or equilibrium. Since equilibria are impossible by assumption, and from above analysis there exist no $\omega$-limit sets in $\bR$, all solutions starting in $\bR$ will converge to $\Omega_{\leq\widehat{M}_\ell}(U_\ell)$.
	
	\textsc{3rd part.} Recall that from Corollary \ref{cor:local stability}, the set $\Omega_{\leq\widehat{M}_\ell}(U_\ell)$ is contained in the basin of attraction of the origin, where $\widehat{M}_\ell=\min\{\delta_\ell(M_\ell),N_\ell\}$, $M_\ell$ and $N_\ell$ are defined in Corollary \ref{cor:local stability}, and $\gamma_\ell$ is given by Assumption \ref{hyp:Lyapunov inequalities}.

\textsc{4th part.} From the above discussion, the origin is locally stable and globally attractive for \eqref{eq:general system}. Thus, it is globally asymptotically stable for \eqref{eq:general system}. This concludes the proof.
\end{proof}

\begin{remark}
	Note that, if $\Omega_{\leq\widehat{M}_\ell}(U_\ell)=\{0\}$, then $\bR=\Omega_{\leq\widetilde{M}_g}(U_g)$ is a simply connected region, and the second part of the proof of Theorem \ref{thm:global stability} can be reduced to the proof of the known Bendixson's criterion.
\end{remark}

\subsection{Proof of Claim \ref{claim:existence of sub level set and inclusion}}\label{sec:proof of claims}

Let $c$ be a positive real number\footnote{ Such a positive real number always exist. Otherwise, for all $n\geq \bbN$, there exists $y_n$ such that $y_n\in\Omega_{\leq 1/n}(U)$ and $y_n\notin\Omega_{\leq\overline{M}}(V)\times\Omega_{\leq\overline{N}}(W)$.
Since $U$ is proper,  $\Omega_{\leq 1/n}(U)\subset\Omega_{\leq 1}(U)$ is compact.
 Hence, there exists $\{y_{n_j}\}_{j\in\bbN}\subset\{y_n\}_{n\in\bbN}$ such that $y_{n_j}\xrightarrow[]{j\to\infty} y^\ast$ and $U(y^*)=0$. From the positive definiteness of $U$, $y^\ast=0$. Consequently, $y_{n_j}$ is a sequence converging to zero and outside $\Omega_{\leq\overline{M}}(V)\times\Omega_{\leq\overline{N}}(W)$. 
 This is impossible since this set is a neighborhood of the origin.} such that $\Omega_{\leq c}(U)\subset\Omega_{\leq\overline{M}}(V)\times\Omega_{\leq\overline{N}}(W)$.

In the first part, it will be shown that, for all $(x,z)\in\bS$,
	\begin{equation}\label{eq:included}
		\hspace{-1pt}U(x,z)\leq\widehat{M}\Rightarrow \max\{V(x),W(z)\}\leq\min\{\overline{M},\overline{N}\}.\hspace{-4pt}
	\end{equation}
	In the second part, it will be shown that, for all $(x,z)\in\bS$,
	\begin{equation}\label{eq:includes}
		\hspace{-1pt}\widetilde{M}\leq U(x,z)\Rightarrow \max\{\underline{M},\underline{N}\}\leq\min\{V(x),W(z)\}.\hspace{-4pt}
	\end{equation}

\noindent Part 1.
 $U(x,z)\leq\widehat{M}$. This implies $U(x,z)=$ \\ $\max\{\sigma(V(x)),W(z)\}\leq$ $\widehat{M}=\min\{\delta(\overline{M}),\overline{N}\}$.
		
		Assume that $\max\{\sigma(V(x)),W(z)\}=\sigma(V(x))$ and $\min\{\delta(\overline{M}),\overline{N}\}=\delta(\overline{M})$. This implies $\sigma(V(x))\leq\delta(\overline{M})$. From \eqref{eq:proof:set glob att:delta < sigma < gamma inv}, $V(x)\leq\sigma^{-1}\circ\delta(\overline{M})<\overline{M}$. Assume now that $\max\{\sigma(V(x)),W(z)\}=W(z)$ and $\min\{\delta(\overline{M}),\overline{N}\}=\delta(\overline{M})$. This implies $W(z)\leq\delta(\overline{M})\leq\overline{N}$. The other two cases are straightforward. Thus, \eqref{eq:included} holds. Hence, $\Omega_{\leq\widehat{M}}(U)\subset(\Omega_{\leq\overline{M}}(V)\times\Omega_{\leq\overline{N}}(W))$;
		
\noindent Part 2. $\widetilde{M}\leq U(x,z)$. This implies $\widetilde{M}=\max\{\gamma^{-1}(\underline{M}),$ $\underline{N}\}\leq U(x,z)=\max\{\sigma(V(x)),W(z)\}$.
		
		Assume that, $\max\{\gamma^{-1}(\underline{M}),\underline{N}\}=\gamma^{-1}(\underline{M})$ and $\max\{\sigma(V(x)),W(z)\}=\sigma(V(x))$. This implies $\gamma^{-1}(\underline{M})\leq\sigma(V(x))$. From \eqref{eq:proof:set glob att:delta < sigma < gamma inv}, $\underline{M}\leq \gamma\circ\sigma(V(x))<V(x)$. Assume now that, $\max\{\gamma^{-1}(\underline{M}),\underline{N}\}=\gamma^{-1}(\underline{M})$ and $\max\{\sigma(V(x)),W(z)\}=W(z)$. This implies $\underline{N}\leq \gamma^{-1}(\underline{M})\leq W(z)$. The other two cases are straightforward. Thus, \eqref{eq:includes} holds. Hence, $(\Omega_{\leq\underline{M}}(V)\times\Omega_{\leq\underline{N}}(W))\subset\Omega_{\leq\widetilde{M}}(U)$;
	
	Since \eqref{eq:max < min} is a strict inequality, from the continuity and surjectivity of $U$, there exists $(x,z)\in\bS$ such that $\widetilde{M}\leq U(x,z)\leq\widehat{M}$. From \eqref{eq:included} and \eqref{eq:includes}, $\widetilde{M}\leq U(x,z)\leq\widehat{M}\Rightarrow \max\{\underline{M},\underline{N}\}\leq\min\{V(x),W(z)\}\leq\max\{V(x),W(z)\}\leq\min\{\overline{M},\overline{N}\}$. Thus, the inclusion \eqref{eq:inclusion} holds. This concludes the proof.\hfill$\blacksquare$

\section{Conclusion}\label{sec:conclusion}

Systems for which the small gain theorem cannot be used, a sufficient condition for the stability of the resulting interconnected system is proposed. The approach consists in verifying if the small gain condition holds in two different regions of the state space: a local and a non-local. In the gap between both regions, assuming mild properties on the vector field, a sufficient condition ensuring the convergence of solutions, for almost every initial condition, is provided. An approach is proposed for planar system for which Bendixson's criterion does not hold. Two examples illustrate the results.

 The authors plan to extend the proposed approach for the case in which, in a countable number of intervals, the small gain condition holds and, between such intervals, a condition ensuring the absence of $\omega$-limit set holds.

{\bfseries Acknowledgements.} The authors thank the anonymous reviewers for suggestions and fruitful discussions.

\bibliographystyle{plain}
\bibliography{Library}

\end{document}